\documentclass[10pt,journal]{IEEEtran}
\usepackage{graphicx}      
\usepackage{amsmath}           
\usepackage{amssymb}
\usepackage{algorithm}
\usepackage{algpseudocode}
\usepackage{hyperref}
\usepackage{color}
\usepackage{cite}
\usepackage{enumerate}

\usepackage[latin1]{inputenc}  

\usepackage{amsmath}
\usepackage{amssymb}
\usepackage{epstopdf}
\usepackage{subcaption}
\usepackage{cite}

\usepackage{amsthm}

\theoremstyle{definition}
\newtheorem{dfn}{Definition}
\theoremstyle{remark}
\theoremstyle{plain}

\theoremstyle{remark}

\DeclareMathOperator{\E}{E}

\DeclareMathOperator*{\argmin}{arg\,min}

\DeclareMathOperator{\col}{col}
\newcommand{\mbf}[1]{\mathbf{#1}}
\newcommand{\mbs}[1]{\boldsymbol{#1}}
\newcommand{\what}[1]{\widehat{#1}}
\newcommand{\wtilde}[1]{\widetilde{#1}}

\newcommand{\Ehat}{\what{\E}}
\newcommand{\T}{\top}
\newcommand{\out}{y}
\newcommand{\cv}{\mathbf{x}}
\newcommand{\group}{z}
\newcommand{\data}{\mathcal{D}}
\newcommand{\0}{\mathbf{0}}
\newcommand{\1}{\mathbf{1}}

\newcommand{\msub}{\mu}

\newcommand{\reg}{\mbs{\phi}}
\newcommand{\regm}{\mbs{\psi}}
\newcommand{\weight}{\mbf{w}}

\newcommand{\cvcov}{\mbs{\Sigma}}
\newcommand{\cdf}{F}
\newcommand{\groupc}{\tilde{z}}

\newcommand{\resid}{r}
\newcommand{\score}{\pi}
\newcommand{\outconf}{\wtilde{\out}}
\newcommand{\ygrid}{\wtilde{\mathcal{Y}}}

\begin{document}

\title{Model-Robust Counterfactual Prediction Method}
\author{Dave Zachariah and Petre Stoica\thanks{Both authors are at
    Uppsala University. E-mail: \texttt{dave.zachariah@it.uu.se}. This work has been partly supported by the Swedish Research
Council (VR) under contracts 621-2014-5874 and 2016-06079}}

\maketitle

\begin{abstract}
We develop a novel method for counterfactual analysis based on observational data using prediction
intervals for units under different exposures. Unlike methods
that target heterogeneous or conditional average treatment effects of
an exposure, the proposed approach aims to take into account the irreducible
dispersions of counterfactual outcomes so as to quantify
the relative impact of different exposures.
The prediction intervals are
constructed in a distribution-free and model-robust manner based on
the conformal prediction approach. The computational obstacles to this
approach are circumvented by leveraging properties of a tuning-free
method that learns sparse additive predictor models for
counterfactual outcomes. The method is illustrated using both real and synthetic data.
\end{abstract}


\section{Introduction}

In many casual inference problems, the unit of analysis
is subject to an exposure, indexed by $\group$, and is associated with a
continuous outcome (or response) $\out$. For instance, an exposure $\group
\in \{ 0, 1 \}$ may correspond to `not receiving' or `receiving' medication. The
inferential question is then typically posed in counterfactual terms:
\begin{quote}$(\ast)$ ``What would the outcome have been, had the unit 
  been assigned to a different exposure $\wtilde{\group} \neq \group$?''
\end{quote}
The ability to address this question using observational data is
relevant in a wide variety of fields, including clinical trials,
epidemiology, econometrics, policy evaluation, etc. \cite{Hernan&Robins2018_causalinf}

\begin{figure}[!h]
\centering
   \begin{subfigure}[b]{0.46\columnwidth}
   \includegraphics[width=1\linewidth]{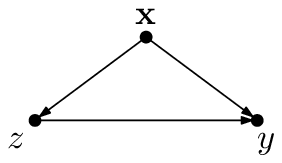}
\end{subfigure}
\begin{subfigure}[b]{0.46\columnwidth}
   \includegraphics[width=1\linewidth]{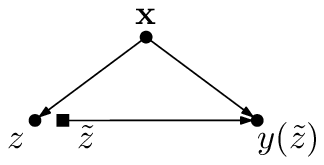}
\end{subfigure}
\caption[TEST]{Dependency graphs, where $\bullet$ and {\tiny
    $\blacksquare$} represent random and deterministic variables,
  respectively. Left: Observed outcome from data-generating
  process. Right: Counterfactual outcome when assigning exposure
  $\groupc$.}
\label{fig:graphs}
\end{figure}

Each unit is typically associated with a range of covariates
(or features), collected in a vector $\cv$, which may affect its
outcome and/or exposure assignment. When $\cv$ contains all variables that
simultaneously affect both $\out$ and $\group$, it is possible to provide
 causal interpretations from observed data. The onus is on the
researcher to include such potentially confounding variables
\cite{Morgan&Winship2014_counterfactuals}. Under this standard condition, the
dependencies between exposure, outcome and covariates can be encoded by a
graph as
in Figure~\ref{fig:graphs} along with an associated joint distribution $p(\cv, \out, \group)$.

The counterfactual question $(\ast)$ can now be directly translated into
predicting the counterfactual (or potential) outcome,
$y(\groupc)$, had the unit been set to $\groupc$,
thus overriding the covariate-dependency of the exposure assignment
\cite{Neyman1923_applications,Rubin1974_estimating,Pearl2009_causality,Pearl2009_causalinstatistics,Sjolander2012_language,Greenland2012_causal,PetersEtAl2017_elements}. The
resulting dependencies for $\out(\groupc)$ can also be encoded in
a graph shown in Figure~\ref{fig:graphs} with an associated joint distribution
 $p( \cv, \out(\tilde{\group}), \group )$
 \cite{Richardson&Robins2013_swigsprimer,Hernan&Robins2018_causalinf}. 

The counterfactual outcome $\out(\tilde{\group})$ is a random variable
and the targeted quantity in most prior works has been the difference
between average outcomes for exposures $\groupc=0$ and $1$, i.e.,
\begin{equation*}
\tau = \E\big[ \: \out(1) \: \big] - \E\big[ \: \out(0) \:\big],
\end{equation*}
which averages out $\cv$
\cite{Morgan&Winship2014_counterfactuals,Imbens2004_nonparametric}.
Using observational data from $n$ units,
$$\data =\bigl \{ (\cv_1, \out_1, \group_1), \; \dots, \;  (\cv_n,
\out_n, \group_n) \bigr \},$$ 
the target $\tau$ is identifiable assuming that the units are drawn independently from the data
generating process $p(\cv, \out, \group)$ and that there is an overlap
of covariates for all exposure types $p(\group | \cv) > 0$ \cite{Wasserman2004_allofstats,Imbens&Rubin2015_causal}. Many methods that estimate $\tau$, model either the
outcome of each exposure type or the exposure selection mechanism as functions
of $\cv$. A central inferential task is to provide confidence
intervals (CI) for the estimate $\what{\tau}$. Much effort has been
made to formulate model-robust methods for this task
as well as extending them to the case of high-dimensional $\cv$ so as
to include a large number of potential confounders, cf. \cite{Robins&Ritov1997_doubly,BelloniEtAl2014_inference,Farrell2015_doublyrobusthighdim,ChernozhukovEtAl2017_debiasedtreatment}. 

For the counterfactual question $(\ast)$, it is however more relevant to
compare the covariate-specific outcomes directly, rather than
averaging them over $\cv$,
cf. \cite{Rothwell1995_clinicaltrial,Kent&Hayward2007_limitations,WeissEtAl2015_ite}. Consequently,
the focus of recent methods has been the covariate-specific effect,
\begin{equation*}
\begin{split}
\tau(\cv) &= \E\big[ \: \out(1) \: | \: \cv \: \big] - \E\big[ \:
\out(0) \: | \: \cv \: \big] \\
&= \msub_{1}(\cv) - \msub_{0}(\cv)
\end{split}
\end{equation*}
also referred to as the `conditional average treatment effect'. Since $p(\out(\tilde{\group}) | \cv) = p(\out  | \cv, \group =
\tilde{\group})$ follows from the dependency structure, it is possible to learn flexible regression models of $\msub_{\tilde{\group}}(\cv)$ using a subset of the observational data,
\begin{equation*}
\data_{\tilde{\group}} =\bigl \{ (\cv_i, \out_i ) \bigr \}, \quad
\text{where} \quad (\cv_i, \out_i) \sim p(\:\cv, \out \: | \: \group
 = \tilde{\group}\:).
\end{equation*} 
The average effect is then estimated as
$\what{\tau}(\cv) = \what{\msub}_1(\cv) - \what{\msub}_0(\cv)$. Using tree-based models it is possible to derive CIs for $\what{\tau}(\cv)$ that are asymptotically valid, cf. \cite{HastieEtAl2013_elements,Hill2011_bart, Wager&Athey2017_estimation, KunzelEtAl2017_meta}. 

A fundamental limitation of targeting $\tau(\cv)$ is that the irreducible dispersions of the
counterfactual outcomes $\out(1)$ and $\out(0)$ are omitted. While correctly
 inferring that, say, $\msub_{1}(\cv) > \msub_{0}(\cv)$, it may still be the case that $\out(0)$
 frequently exceeds the value of $\out(1)$. Such
 considerations are important in applications where different expoures
 involve differential risks. Then merely reporting $\what{\tau}(\cv)$
 and a CI omits valuable information about $\out(1)$ and $\out(0)$.

\begin{figure}
  \begin{center}
    \includegraphics[width=0.53\columnwidth]{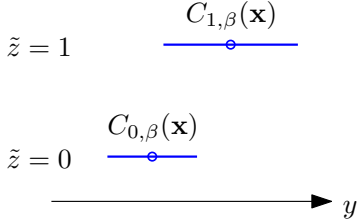}
  \end{center}
  \caption{Counterfactual analysis of units with covariates $\cv$
    assigned to different exposures $\groupc$. Outcome predictions
    $\what{\msub}_0(\cv)$ and $\what{\msub}_1(\cv)$ are denoted by
    circles. The  prediction intervals (PIs) $C_{0, \beta}(\cv)$ and
    $C_{1, \beta}(\cv)$ (lines) incorporate the
    outcomes with a coverage of $100\% \beta$. Together the PIs
    enable an intuitive assessment of the counterfactual
    outcomes. Below we propose using the PIs to quantify how separable the outcomes are under different exposures.}
  \label{fig:space_observable}
\end{figure}

 In  this paper, we aim to address this limitation by considering
 $\what{\msub}_{\tilde{\group}}(\cv)$ as a predictor of $\out(\groupc)$ for a unit with covariates $\cv$. Then using
 prediction intervals (PIs) for the predictors
 \cite{Geisser1993_predictive}, the irreducible dispersions of the
 outcomes can be taken into account in the counterfactual analysis,
 cf. Fig.~\ref{fig:space_observable}.
A major challenge is obtaining PIs that are valid when the
data-generating process is unknown.  
Here we consider the general conformal approach
 of \cite{VovkEtAl2005_algorithmic,BalasubramanianEtAl2014_conformal}
to construct PIs with valid coverage properties in
 distribution-free manner. Each point in the interval is then
 constructed by re-fitting the predictor for the corresponding
 exposure group. This becomes  computationally prohibitive using
 complex tree-based models or fitting methods that require parameter
 tuning, especially when $\cv$ is high-dimensional. This obstacle is
 circumvented using the method proposed below.

Our contributions in this paper are the following. We propose the use
of prediction intervals for covariate-specific counterfactual analysis
of each exposure type and define a measure of their relative impact
that provides additional information to merely comparing $\what{\msub}_0(\cv)$ and $\what{\msub}_1(\cv)$. We then learn sparse predictor models $\what{\msub}_{\tilde{\group}}(\cv)$, with corresponding PIs $C_{\groupc}(\cv)$, that automatically adapt to nonlinearities by leveraging the computational properties of the \textsc{Spice} predictor approach \cite{ZachariahEtAl2017_online}. This obviates the need for cross-validation or other tuning techniques. Since the conformal PIs also exhibit marginal coverage properties, even when lacking a correctly specified model of the data-generating process, the resulting method provides a model-robust means of counterfactual analysis.










\emph{Notation:} $\| \cdot \|_1$ and $\odot$ denote the
$\ell_1$-norm and Hadamard product. The cardinality of a set $\mathcal{D}$
is $|\mathcal{D}|$. The operator $\col\{ \cv_1, \dots, \cv_k  \}$
stacks all elements into a single column vector. $\Ehat[ f(\cv,\out)
|\groupc] = |\data_{\groupc}|^{-1} \sum_{i} f(\cv_i, \out_i)$ is the
empirical mean of $f(\cv,\out)$ over all pairs $(\cv, \out)$ in $\data_{\groupc}$.

\emph{Remark:} Code for the proposed method is available at \url{https://github.com/dzachariah/counterfactual}.

\section{Counterfactual prediction method}
\label{sec:method}

A predictive approach to counterfactual analysis is readily
generalized to multiple exposure types $\group \in \{0 , \dots, K-1
\}$. We illustrate this as we proceed to define a measure of the
separability of counterfactual outcomes under different exposures.

\subsection{Counterfactual confidence}

\begin{figure*}[!h]
\centering
   \begin{subfigure}[b]{0.32\textwidth}
   \includegraphics[width=1\linewidth]{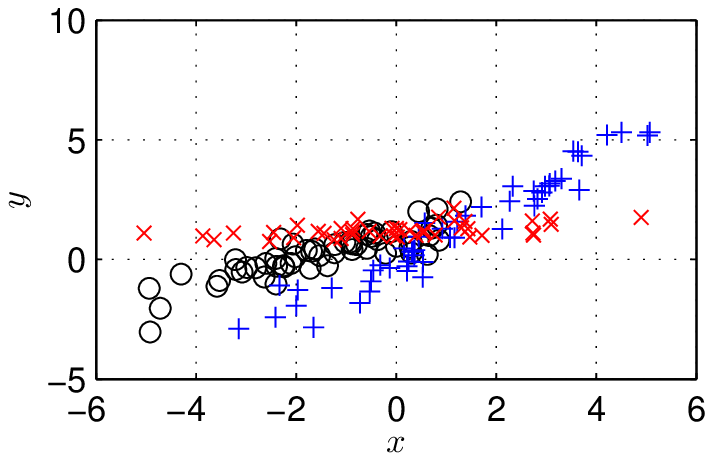}
\end{subfigure}
$\:$
\begin{subfigure}[b]{0.32\textwidth}
   \includegraphics[width=1\linewidth]{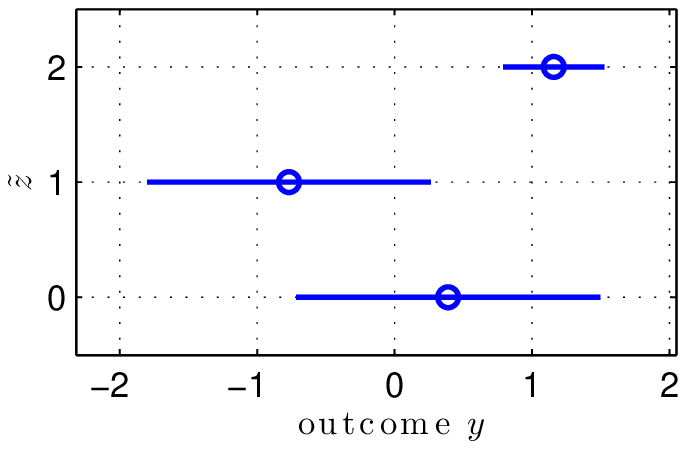}
\end{subfigure}
$\:$
\begin{subfigure}[b]{0.32\textwidth}
   \includegraphics[width=1\linewidth]{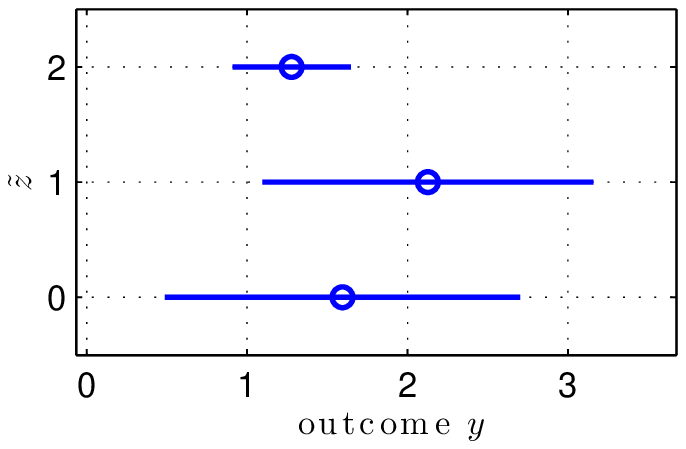}
\end{subfigure}
\caption[TEST]{Counterfactual analysis with $K=3$
  exposure types and scalar covariate $x$, using method proposed on subsequent sections. Left: Datasets $\data_{\groupc}$ for exposures
  $\groupc = 0$ ($\circ$), $1$ ($+$) and $2$ ($\times$). The resulting
  predictions $\what{\msub}_{\groupc}(x)$ along with $90\%$-intervals
  $C_{\groupc,0.90}(x)$. Center: $x=-1$. Right: $x=2$. }
\label{fig:example}
\end{figure*}

Consider an observational study with $K=3$ exposure types and a
scalar covariate $x$, as illustrated in
Fig.~\ref{fig:example}. For a given $x$, let $\what{\msub}_{\groupc}(x)$ and 
prediction intervals $C_{\groupc, \beta}(x)$ be the counterfactual
prediction and PI for each exposure $\groupc \in \{ 0, 1, 2 \}$. We
propose the following measure of the impact of exposures on the outcomes.
\begin{dfn} A covariate-specific comparison between predicted outcomes for exposures $\groupc$ and $\group$ is said to have $100 \beta$\%
\emph{counterfactual confidence}, where $\beta$ is the largest value
for which their PIs do not overlap. That is, the PIs are mutually exclusive
\begin{equation}
C_{\groupc,\beta}(\cv) \cap  C_{\group,\beta}(\cv) =  \emptyset.
\label{eq:defcounterfactual}
\end{equation}
A high confidence asserts that the counterfactual outcomes are highly
separable and thus the impact of the exposures are more
distinctive. This provides additional information to the size of the exposure effect measured as
$\what{\msub}_{\groupc}(\cv) - \what{\msub}_{\group}(\cv)$.
\end{dfn}
In Fig.~\ref{fig:example} we see that
for a unit with covariate $x=-1$, the counterfactual confidence when comparing exposures
$1$ and $2$ is greater than $90\%$, indicating highly separable outcomes. The
confidences for the pairwise are tabulated pairwise below
for covariate $x=-1$ (left) and $x=2$ (right):  
\begin{center}
\begin{tabular}{c|cc}
\hline
$\tilde{\group}$&$0$&$1$\\
\hline
1 & 81\% & -- \\
2 & 64\% & 96\% \\
\end{tabular}
 \qquad
\begin{tabular}{c|cc}
\hline
$\tilde{\group}$&$0$&$1$\\
\hline
1 & 24\% & -- \\
2 & 45\% & 58\% \\
\end{tabular}
\end{center}
It is seen that for $x=-1$, the separability of counterfactual outcomes
can be asserted with greater confidence than for $x=2$. This is
corroborated by comparing the datasets shown in Fig.~\ref{fig:example}.



\subsection{Conformal prediction intervals}

Consider a regression model class
\begin{equation*}
\mathcal{M}_{\msub} = \Big\{ \msub : \msub(\cv) = \reg^\T(\cv)
\weight, \; \weight \in \mathbb{R}^{p+1} \Big\},
\end{equation*}
parameterized by $\weight$. In
subsection~\ref{sec:regressionfunction}, we specify a flexible
regressor vector $\reg(\cv)$ that adapts to nonlinearities in the
data. When the model class is
well-specified it includes the unknown mean function
$\msub_{\group}(\cv) \in \mathcal{M}_{\msub}$.

To learn a model in $\mathcal{M}_{\msub}$ from $\data_\group$, we
build upon the tuning-free \textsc{Spice}-method
\cite{ZachariahEtAl2017_online} where the learned weights are
defined as
\begin{equation} \what{\weight} \: = \: \argmin_{\weight} \: \sqrt{\Ehat\big[ \: | \out -
  \reg^\T(\cv)\weight|^2 \: \big| \: \groupc \: \big]} +\: \| \mbs{\varphi} \odot \weight  \|_1,
\label{eq:coeffmin_spice}
\end{equation}
and the elements of $\mbs{\varphi}$ are given by
\begin{equation*}
\varphi_j = \begin{cases}\sqrt{ \frac{ \Ehat \big[ \: \phi^2_{j}(\cv)
      \:  | \: \groupc \: \big]}{|\data_{\groupc}|} }, & j = 1, \dots,
  p\\
0, & j=0.
\end{cases}
\end{equation*} 
The solution \eqref{eq:coeffmin_spice} can be computed sequentially
for each sample in $\data_{\groupc}$, using a coordinate descent
algorithm with a runtime that scales as $\mathcal{O}(|\data_{\groupc}|
p^2)$. We now leverage this property for the construction of conformal PIs.

The principle behind general conformal prediction can be described as follows
\cite{VovkEtAl2005_algorithmic,BalasubramanianEtAl2014_conformal}: For the covariate of interest $\cv$,
consider a new sample $(\cv, \outconf)$ where $\outconf$ is a variable. Then quantify how well this sample conforms to the
observed data $\data_{\groupc}$ via the learned model
\eqref{eq:coeffmin_spice}. All points $\outconf$ that conform
well, form a prediction interval. The conformity is quantified by including $(\cv, \outconf)$ in
the learned model \eqref{eq:coeffmin_spice}, which is achieved by a
sequential update in
$\mathcal{O}(p^2)$. Then, following \cite{LeiEtAl2017_distribution}, we define a measure
\begin{equation}
\score(\outconf) = \frac{1}{|\data_{\groupc}|+1}\Big( 1+ \sum_i I\big\{
          \resid_i \leq |\outconf - \reg^\T(\cv) \what{\weight}|
          \big\} \Big),
\label{eq:score}
\end{equation}
where $\resid_i = |\out_i - \reg^\T(\cv_i) \what{\weight}|$ is the $i$th
fitted residual and $I\{ \cdot \}$ is the indicator function. The measure is bounded between 0 and 1, where lower
values correspond to higher conformity. We construct $C_{\groupc, \beta}(\cv)$
by varying $\outconf$ over a set of grid points $\ygrid$, as
summarized in 
Algorithm~\ref{alg:conformal}. By leveraging the computational
properties of the learning method, the prediction interval
is computed with a total runtime of
$\mathcal{O}(|\ygrid|p(p + |\data_{\groupc}|))$. The range of $\ygrid$ is set to exceed that of the outcomes in the observed dataset. A point prediction  $\what{\msub}_{\groupc}(\cv)$ is
obtained as the minimizer of $\score(\outconf)$.
\begin{algorithm}
  \caption{: Conformal prediction interval} \label{alg:conformal}
\begin{algorithmic}[1]
    \State Input: covariate $\cv$, target coverage $\beta$ and data $\data_{\groupc}$
    \For{ all $\outconf \in \ygrid$}
        \State Update $\what{\weight}$ using $(\cv, \outconf)$
        \State Compute $\{ \resid_i\}$ and $\score(\outconf)$ in \eqref{eq:score}
    \EndFor
    \State Output: $C_{\groupc, \beta}(\cv) = \big\{ \outconf \in
    \ygrid \: : \:  (n+1)\score(\outconf) \leq \lceil
    \beta(n+1) \rceil \big\}$
\end{algorithmic}
\end{algorithm}

Despite the fact that no dispersion model of the data generating process is required, the
resulting prediction intervals exhibit valid coverage
properties. When the model class $\mathcal{M}_\msub$ is well-specified,
the interval exhibits asymptotic conditional coverage, that is,
\begin{equation*}
\Pr\big\{ \: \out \in C_{\groupc,\beta}(\cv )   \: \bigl|  \: \cv, \groupc   \: \big
\}  \: =  \: \beta + o_P(1), 
\end{equation*}
under certain regularity conditions \cite[thm.~6.2]{LeiEtAl2017_distribution}. More generally, $C_{\groupc,\beta}(\cv )$ is calibrated to
ensure marginal coverage
\cite[thm.~2.1]{LeiEtAl2017_distribution}
\begin{equation*}
\Pr\big\{ \: \out \in C_{\groupc,\beta}(\cv )   \: \big|  \: \groupc   \: \big
\}  \: \geq  \: \beta.
\end{equation*}
Note that this does not
require a well-specified model class $\mathcal{M}_\msub$. In other
words, the more accurate the learned prediction model, the tighter the
prediction interval but its marginal coverage property remains not matter if the model is correct or not. This confers a robustness property to the proposed inference method in cases when $\msub_{\groupc}(\cv) \not\in \mathcal{M}_\msub$.

\subsection{Sparse additive predictor models}
\label{sec:regressionfunction}

To learn accurate prediction models, we now turn to the specification
of the regressor vector in $\mathcal{M}_\msub$. Let $d$ denote the
dimension of $\cv$ and consider the additive model
\begin{equation*}
\reg(\cv) = \col\big\{1, \: \regm_1(x_1), \: \dots,  \: \regm_d(x_d) \big\},
\end{equation*}
which enables interpretable component-wise predictors
\cite{HastieEtAl2013_elements}. When a covariate $x$ is categorical,
we use a standard basis vector $\regm(x) \equiv \mbf{e}_k$ for
category $k \neq 0$ and $\0$ when $k=0$. In the case of binary categories, we
simply have $\regm(x) \equiv 1$.

For noncategorical covariates $x$, such as continuous or count
variables, however, we propose a data-adaptive piecewise linear
model. Let the empirical quantile function be
\begin{equation*}
\what{\cdf}^{-1}( q) = \inf\big\{ x: \what{\cdf}(x) \leq q \big\},
\end{equation*}
where $\what{\cdf}(x)$ is the empirical cumulative density function of
$x$ obtained from the full observational data $\data$. Then we specify
the model for $x$ by the $m$-dimensional vector
\begin{equation*}
[\regm(x)]_k =
\begin{cases}
(x - c_k)_+ ,& k =1, \dots, m-1\\
(x - c_k)_+I\{x \leq c_{k+1}\}  & \\
\:+\: c_{k+1}I\{x > c_{k+1}\}, & k = m
\end{cases}
\end{equation*}
where $m$ is a specified number of knots or breakpoints defined as
$c_k = \what{\cdf}^{-1}( \frac{k-1}{m} )$. The resulting model yields
finer resolution line segments where data density is high and is capable
of capturing continuous nonlinear responses with respect to $x$.

The regularization term in \eqref{eq:coeffmin_spice} leads to sparse
solutions and thus the learning method yields a set of sparse
additive predictor models $\what{\msub}_{\groupc}(\cv) =
\what{\msub}_{\groupc}(x_1) + \dots + \what{\msub}_{\groupc}(x_d)$ with nonlinear responses in a
data-adaptive manner, cf. \cite{RavikumarEtAl2007spam}.\footnote{The
  predictive performance can be related to the best subset predictor, see
\cite{ZachariahEtAl2017_online} for more details.} Moreover, the
method takes into account the amount of data observed within each line
segment for each exposure type, and controls the model complexity accordingly. This mitigates overfitting when there is unequal
training data or covariate imbalances across exposure types,
cf. \cite{KunzelEtAl2017_meta}.

The integer $m$ determines the maximum resolution of the nonlinear
model. A high $m$ enables higher predictive accuracy provided enough
data is available to learn the nonlinear response. For instance,
suppose $n$ denotes the size of the smallest dataset $\data_{\groupc}$ and
$d'$ and $d''$ denote the number of binary and continuous covariates,
respectively. Then $p = d' +  md''$ and a natural upper limit is $m
\leq \max( \lceil\frac{n - d'}{d''} \rfloor, 1)$.

\section{Numerical experiments}
\label{sec:results}

In this section we demonstrate the proposed counterfactual prediction
approach by means of three examples. In the following examples, we
consider $K=2$ exposure types.

\subsection{Nonlinear effects}

\begin{figure*}[!h]
\centering
   \begin{subfigure}[b]{0.49\textwidth}
   \includegraphics[width=1\linewidth]{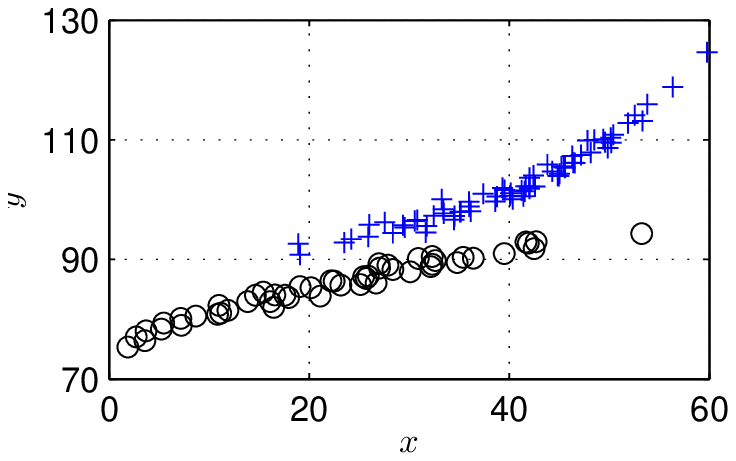}
\end{subfigure}
$\:$
\begin{subfigure}[b]{0.49\textwidth}
   \includegraphics[width=1\linewidth]{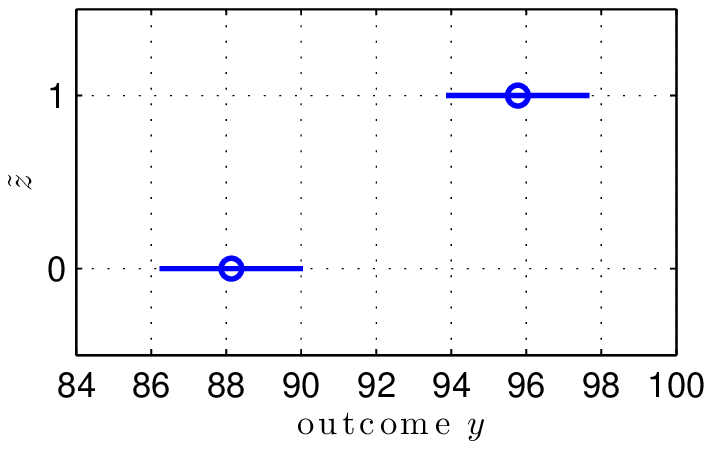}
\end{subfigure}
\caption[TEST]{Left: Dataset $\data$ with $n=120$ samples. Exposures
  $\groupc = 0$ ($\circ$) and $1$ ($+$). Right: Predictions with $90\%$-prediction intervals at $x = 30$.}
\label{fig:nonlinear}
\end{figure*}

We consider the example in
\cite{Hill2011_bart}, which allows for a comparison with both
tree-based and linear models. For each
unit, the exposure $\group$ is assigned with equal probability. Then the
covariate $x$ (with $d=1$) is  drawn as
$$x|(\group = 0) \sim \mathcal{N}(40,10^2) \quad \text{and} \quad x|(\group = 1) \sim \mathcal{N}(20,10^2)$$
and the counterfactual outcomes as
\begin{equation}
\begin{split}
y(0)|x &\sim \mathcal{N}(72+3\sqrt{|x|},1) \quad \text{and}\\
y(1)|x &\sim \mathcal{N}(90+\exp(0.06x),1).
\label{eq:model_nonlinear}
\end{split}
\end{equation}
A simulated observational dataset $\data$ with $n=120$ is
illustrated in Fig.~\ref{fig:nonlinear}. To obtain the predictions in
the figure, we use $m=10$.  For a unit with covariate $x=30$, as an example, we note
that $\what{\msub}_1(30)$ is larger than
$\what{\msub}_0(30)$ and that both confidence
intervals are tight, as is expected by inspecting the data generating process
\eqref{eq:model_nonlinear} at the given covariate. In addition, the
counterfactual confidence is found to be greater than 90\%, indicating
a highly separable counterfactual outcomes as expected.

To illustrate the robustness property of the prediction intervals,
we repeat the experiment using 1000 Monte Carlo simulations. For each
simulation, we generate new data $\data$ and also draw a new unit from
both exposure
groups. For a unit with exposure $\group = 0$, the outcome
is found to belong to interval $C_{0,\beta}(x)$ with probability
0.918 when $\beta = 0.90$.  Similarly, for exposure $\group = 1$, the
outcome is contained in the interval $C_{1,\beta}(x)$ with
probability 0.907. This coverage property holds eventhough the mean
function does not belong to $\mathcal{M}_\msub$.

\subsection{High-dimensional covariates}

The desire to include all potential confounders in the covariate vector $\cv$,
may lead in many applications to dimensions $d$ that can be larger
than $n$ \cite{Farrell2015_doublyrobusthighdim}. To address this
issue, we simulate an experimental setting with $d=200$ covariates but
only $n =100$ samples. Naturally, we have $m=1$. The exposures $\group=0$ and $\group=1$ are assigned with probabilities 0.6 and 0.4, respectively. The
covariates are drawn as
$$\cv|(\group = 0) \sim \mathcal{N}(\0, \cvcov_0 ) \quad \text{and} \quad \cv|(\group = 1) \sim \mathcal{N}(\0,\cvcov_1),$$
where $\cvcov_0$ and $\cvcov_1$ are randomly generated covariance
matrices with unit trace. The matrices have numerical rank $150$ and
are constructed using outer products of Gaussian vectors. This generates highly correlated covariates, as is typical in real applications. The counterfactual outcomes are generated as
\begin{equation}
\begin{split}
y(0)|\cv &\sim \mathcal{N}(x_1 + 5x_{10} + 5x_{20}+0.5, \; 0.5^2) \quad \text{and} \\
y(1)|\cv &\sim \mathcal{N}(x_1 +  x_{10} - x_{30}, \; 0.5^2).
\end{split}
\label{eq:model_highdim}
\end{equation}
However, this is not a problem for the learning method which automatically prunes away irrelevant covariates due to the adaptive regularization in \eqref{eq:coeffmin_spice}.

\begin{figure*}[!h]
\centering
   \begin{subfigure}[b]{0.49\textwidth}
   \includegraphics[width=1\linewidth]{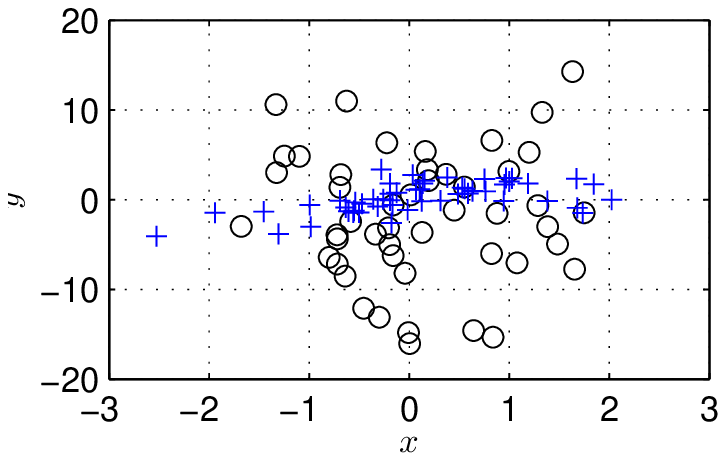}
\end{subfigure}
$\:$
\begin{subfigure}[b]{0.49\textwidth}
   \includegraphics[width=1\linewidth]{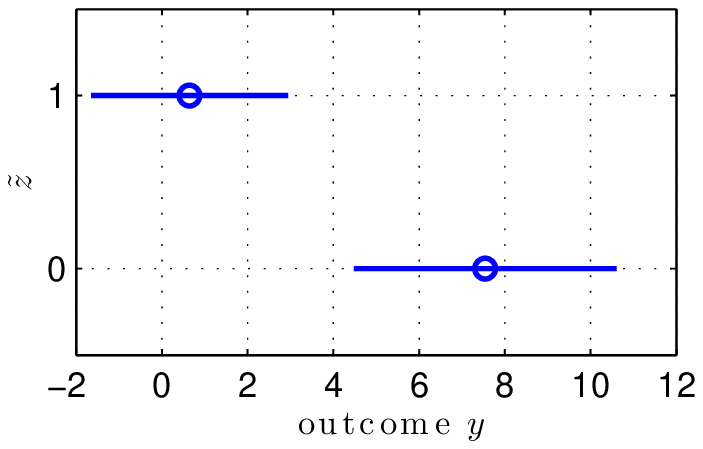}
\end{subfigure}
\caption[TEST]{Left: Dataset $\data$ with $n=100$ samples. Only one
  covariate, $x_1$, is shown for the sake of two-dimensional visualization. Exposures
  $\groupc = 0$ ($\circ$) and $1$ ($+$). Right: Predictions with $90\%$-prediction intervals at $\cv = \1$.}
\label{fig:highdim}
\end{figure*}

A simulated observational dataset $\data$ is shown in Fig.~\ref{fig:highdim}. We also illustrate the predicted outcomes for a unit with all covariates equal to one, $\cv = \1$. We observe that $\what{\msub}_0(\1)$ is considerably larger than
$\what{\msub}_1(\1)$, also when taking into account the prediction
intervals. This is consistent with the data
generating process \eqref{eq:model_highdim} evaluated at the fixed $\cv$. The interval for exposure
$\group = 0$ is also seen to be significantly wider than that for
exposure $\group =1$, reflecting the larger uncertainty of the predicted
outcome. In this case it is possible to assert counterfactual
confidence greater than 90\%.

We repeat this experiment as well to validate the coverage properties of the intervals, using 1000 Monte Carlo simulations. For each
simulation, we generate new data $\data$ and also draw a new unit from
both exposure
groups. For a unit with exposure $\group = 0$, the outcome
is found to be contained in interval $C_{0,\beta}(\cv)$ with probability
0.921 when $\beta = 0.90$.  Similarly, for exposure $\group = 1$, the
outcome is contained in the interval $C_{1,\beta}(\cv)$ with
probability 0.915.

\subsection{Schooling data}

Following the example in \cite{Esarey2017_causal}, we assess the 
effect of schooling on income for adults in the US born in the 1930s, using data from 
\cite{Angrist&Krueger1991_school}. The observed outcome $y$ is the
weekly earnings (on a logarithmic scale) of a subject in 1970. Each
individual is subject to one of two exposures:
$\group = 1$ corresponds to receiving 12 years
of schooling or more and $\group = 0$ corresponds to receiving less
than 12 years. We consider 26 binary covariates in $\cv$. Ten
covariates indicate the year of birth 1930-1939, and
eight indicate the census region. In addition, eight indicators
represent ethnic identification, marital status and whether or not the
subject lives in the central city of a metropolitan area. The observational study consists of $n=329~509$ samples. (See
\cite{Esarey2017_causal} for details.)

\begin{figure*}[!h]
\centering
   \begin{subfigure}[b]{0.32\textwidth}
   \includegraphics[width=1\linewidth]{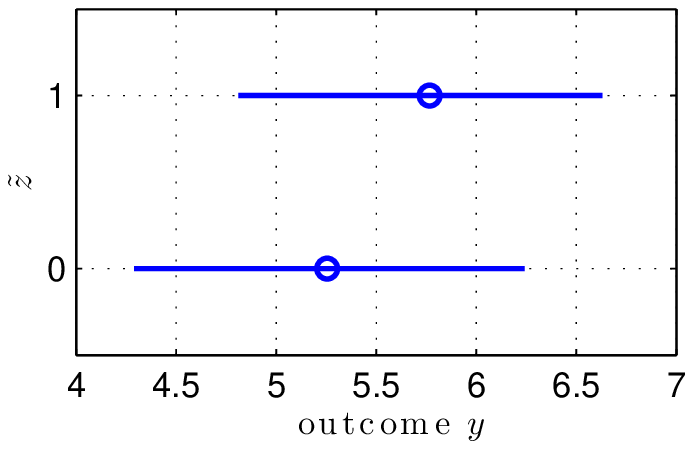}
\end{subfigure}
$\:$
\begin{subfigure}[b]{0.32\textwidth}
   \includegraphics[width=1\linewidth]{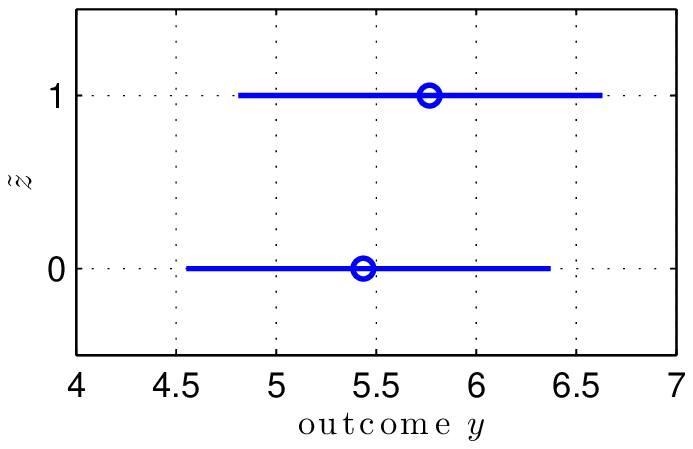}
\end{subfigure}
$\:$
\begin{subfigure}[b]{0.32\textwidth}
   \includegraphics[width=1\linewidth]{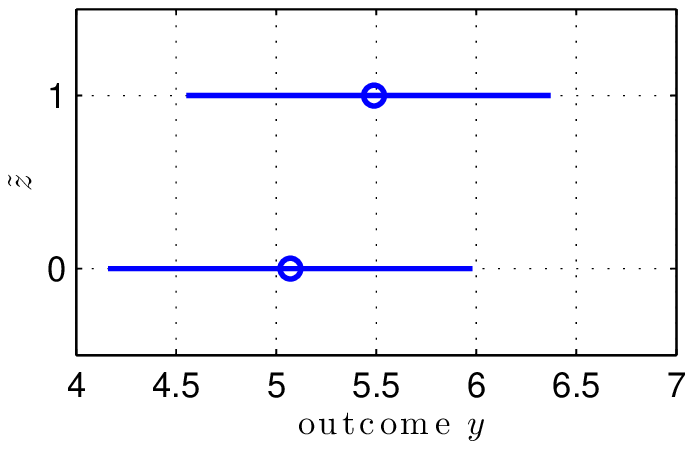}
\end{subfigure}
\caption[TEST]{Counterfactual analysis of schooling ($\groupc$) and
  weekly earnings ($\out$, in logarithmic units) using $90\%$-prediction intervals. Left: Caucasian, unmarried and not in major
  city. Center: Caucasian, married and in city. Right: African-american, married and in city.}
\label{fig:school}
\end{figure*}

Discrete
covariates can be partitioned into separate subgroups, and a
direct inference approach would be to estimate the average outcomes of
exposures 0 and 1 for each group. However, the number of subgroups grows
quickly and there are not sufficient samples in the dataset $\data$ for each subgroup and exposure. Therefore we apply the
proposed method. The predicted outcomes 
are illustrated for subjects in  different covariate groups in
Fig~\ref{fig:school}. All subjects in these subgroups were born in the
same year and came from the same region. The prediction
interval widths are likely to be affected by the very coarse division of schooling
used here, since $\group=0$ includes 0 to 11
years of schooling, which is a substantial variation, while $\group=1$ includes 12 years and more.

The three subgroups are $\cv_1$: Caucasian, unmarried and not in a major
city, $\cv_2$: Caucasian, married and in a major city, and $\cv_3$:
African-American, married, and in a major city. Given that the units are
logarithmic, the differences of predicted earnings,
$\what{\msub}_1(\cv) - \what{\msub}_0(\cv)$, correspond to +52\%,
+26\% and +39\% of weekly earnings, for  $\cv_1$,  $\cv_2$ and $\cv_3$,
respectively. This means that the inferred effect of schooling is greatest
for  $\cv_1$ while considerably less for  $\cv_2$. The prediction
intervals in Fig.~\ref{fig:school} suggest, however, that there is a
considerable dispersion of the outcomes. The
predicted outcome of schooling has a
counterfactual confidence of 33\%, 20\% and 25\% for the three subgroups.
Thus for subgroup $\cv_2$, schooling not only exhibits the lowest
predicted gains but also the least separable counterfactual outcomes
can be asserted. The opposite is true for subgroup $\cv_1$.

The findings appear to be consistent with features of
US society in the 1970s: a Caucasian person in a major
city with a family was expected to have greater access to economic opportunities, such that schooling
experience mattered less to earnings. For the unmarried counterpart
who lived outside of the major city, such alternative opportunities were
fewer so that schooling could have a more significant impact. An
African-American person in a major city with a family represents an
intermediate case.

\section{Conclusions}

We have developed a new method for counterfactual analysis using
observational data based on prediction intervals. The intervals were
used to define a measure of relative separability of counterfactual
outcomes under different exposures. This takes into account the dispersions of
the outcomes and provides additional information to
the difference between predictions.

The intervals were constructed in a
distribution-free and model-robust manner based on the general conformal
prediction approach. The computational obstacles of this approach were
circumvented by leveraging properties of a tuning-free
method that learns sparse additive predictor models for
counterfactual outcomes. We demonstrated the method
using both real and synthetic data.


\bibliographystyle{ieeetr}
\bibliography{refs_causalinf}

\end{document}